\documentclass[12pt]{Tran-l}

 \newtheorem{thm}{Theorem}[section]
 \newtheorem{cor}[thm]{Corollary}
 \newtheorem{lem}[thm]{Lemma}
 \newtheorem{prop}[thm]{Proposition}
 \theoremstyle{definition}
 
 \theoremstyle{remark}
 
 \theoremstyle{definition}

 \newcommand{\CC}{\mathbb{C}}

 \newcommand{\PP}{\mathbb{P}}

\begin{document}

\title{Existence results for rational normal curves}

\author{E. Carlini and M. V. Catalisano}





\begin{abstract}
In this paper we study existence and uniqueness of rational normal
curves in $\PP^n$ passing through $p$ points and intersecting $l$
codimension two linear spaces in $n-1$ points each. If $p+l=n+3$
and the points and the linear spaces are generic, one expects the
curve to exist, but this is not always the case. Our main result
precisely describes in which cases the curve exists and in which
it does not exist.
\end{abstract}

\maketitle

\section{Introduction}

In this paper we study existence and uniqueness of rational normal
curves in $\PP^n$ passing through a given set of points and
intersecting some codimension two linear spaces in a very natural
way. More precisely, we require the curve to intersect each linear
space in $n-1$ distinct points. In this case, we say that the curve
and the linear space are mutually $(n-1)$-secant. We work over the
field of complex number $\mathbb{C}$ and we consider a {\it rational
normal curve} (briefly a {\it rnc}) as a linearly normal embedding of
$\PP^1$.

Our interest in the subject arises from the deeply
intertwined problems of the postulation of schemes and of the dimension
of higher secant varieties, e.g. see the original work of
Terracini \cite{Terracini} and also \cite {Ge}. We give a first
application of our results in this direction in Section
\ref{postulationAPP} where we easily obtain a well known result
by {\AA}dlandsvik \cite{AAdlandsvik} about Segre-Veronese varieties,
 recently re-proposed by
Abrescia in \cite{Abrescia}.

The idea of using rational curves in the study of linear systems
and higher secant varieties is classical. Its importance has been
stressed again in the case of double points schemes and higher
secant varieties to Veronese varieties. To explain this, let
$X\subset\PP^n$ be a double point scheme supported on $p$ generic
points $P_i$'s, i.e. a scheme with defining ideal
\[
I_X=\left(I_{P_1}\right)^2\cap\ldots\cap \left(I_{P_p}\right)^2
\]
where the $I_{P_i}$'s are the ideals of the $P_i$'s. Then, one
wants to determine the Hilbert function of $X$ in some degree $d$,
say $H(X,d)$. There is an expected value for the Hilbert function
 determined by a naive count of conditions, which we will call $h(n,p,d)$. This value is such
that $H(X,d)\leq h(n,p,d)$ and if the points are generic one {\it
expects} equality to hold. In a series of papers, Alexander and
Hirschowitz  determined exactly when equality holds, see
\cite{AH95} and \cite{Chandler}. More precisely, $H(X,d)= h(n,p,d)$ in
all but the following cases:
\begin{itemize}
\item $d=2$;
\item $d=4, (n,p)=(2,5),(3,9),(4,14)$;
\item $d=3, (n,p)=(4,7)$.
\end{itemize}
The $d=4$ cases are easily explained by the existence of quadric
hypersurfaces passing through $5,9$ and $14$ points in
$\PP^2,\PP^3$ and $\PP^4$,
 respectively.
The $d=3$ case requires a subtler explanation involving rational
normal curves. Given a scheme $X\subset\PP^4$ consisting of seven
double points, we do not expect a cubic threefold singular at all
the points to exist, i.e. $h(4,3,7)={4+3\choose 3}=35$. But there
is a rnc $\mathcal{C}$ passing through the points (see Theorem
\ref{castelnuovo}) and the variety of secant lines to
$\mathcal{C}$ is a cubic threefold singular along the curve. Thus,
$H(X,3)<35$ and the variety of secant $\PP^6$'s to the 3-uple
embedding of $\PP^3$ does not have the expected dimension. For a
more detailed account see \cite{Ge}, the introduction of
\cite{RS00} and \cite{Ci01}.

Our research interest was inspired by the following classical
result

\vskip .25cm

\noindent{\bf Theorem \ref{castelnuovo}.}{\it Given
$n+3$ points in $\PP^n$ in generic position, there exists a unique rational
normal curve passing through them.}

\vskip .25cm

This theorem was well known in the late $19^{th}$ century, e.g. it
can be found in works by Bordiga \cite{Bordiga} and Castelnuovo
\cite{Castelnuovo} where it is attributed to Veronese
\cite{Veronese}. We want to mention that even then there were
attempts to generalize  this result, but always in a constructive
way. In the sense that the final goal of those attempts was the
synthetic construction of a curve satisfying certain properties.
Recently Graber and Ranestad in \cite{GraberRanestad}, following
Kapranov \cite{Kapranov}, generalized Theorem \ref{castelnuovo}
 to d-uple Veronese surfaces and applied their results to the study of line arrangements.

In Section \ref{main} we provide a classical proof of Theorem
\ref{castelnuovo}. We note that it is natural to expect that $n+3$
points determine a finite number of rnc's. In fact, the parameter
space of rnc's in $\PP^n$ has dimension $(n+3)(n-1)$. Moreover, a
simple argument on the defining matrix of a rnc shows that the
family of rnc's passing through a given point has dimension
$(n+2)(n-1)$. Thus a point imposes $n-1$ conditions to rnc's.
Theorem \ref{castelnuovo} is nothing more than the proof that
$n+3$ points impose independent conditions.

We can push this kind of argument even further. Take a codimension
two linear space $\Lambda$ in $\PP^n$ and consider all rnc's
intersecting it: this is one condition for rnc's. Thus, if we
consider rnc's intersecting $\Lambda$ in $n-1$ points this
incidence condition imposes again $n-1$ conditions. In conclusion,
passing through a fixed point or intersecting a given codimension
two linear space in $n-1$ not fixed points imposes the same number
of conditions to rnc's. Hence, it is natural to look for
generalizations of Theorem \ref{castelnuovo} involving points and
codimension two linear spaces.

As a first step in this direction, one tries to generalize Theorem
\ref{castelnuovo} in $\PP^3$. The case of twisted cubic curves was
studied in details in the early $20^{th}$ century and even in this simple case it
is clear that the theorem does not generalize in a  straightforward way.
In fact, for a generic choice of four points and two lines,
 via the  count of conditions above, we expect to find at least
one rnc passing through the points and having the lines as chords. But
this is not the case and such a curve does not exist (see Section
\ref{twisted} for more details).

The classic approaches to the $\PP^3$ case use extremely {ad hoc}
arguments which do not generalize easily to higher dimension. We
develop a more general framework where the twisted cubic case and
the general situation can  both be studied. As a result of our
analysis we obtain the following

\medskip
\noindent {\bf Theorem \ref{final}.} {\it Let $n,p$ and $l$ be positive
integers such that
$$n\geq
3, \ \ p\geq 1 \mbox{ and} \ \ p+l=n+3 .$$
 Choose  $p$ points in $\PP^n$ and $l$
codimension two linear spaces in generic position. Then, only for
the values
\[
(p,l)=(n+3,0),(n+2,1),(3,n),(2,n+1),(1,n+2)
\]
does there
exist a unique rational normal curve passing through the points and
$(n-1)$-secant to the linear spaces. In the other cases, that is for
 $p\geq 4$ and $l\geq 2$, no such curve exists.
}

\medskip

 The $(p,l)=(n+3,0)$ case is just
Theorem \ref{castelnuovo}, while the $(p,l)=(n+2,1)$ case
 and the non-existence results
are, as far as we know, original (see Propositions
\ref{onePOINTlessPROP} and \ref{donotexist}). The result for
$(p,l)=(3,n)$ is just Steiner's construction, see e.g. \cite[pg.
528]{GH}, but we provide a different proof not using the classical
construction. The cases $(p,l)=(2,n+1),(1,n+2)$ were studied by
Todd in  \cite{Todd} and by Veneroni in \cite{Veneroni}. Todd
provides a proof of both the results for $n=4$ which we briefly
sketch in the proofs of Propositions \ref{3points} and
\ref{2points}. In \cite{Todd}, the author also claims that the
results extend for any $n$, but no proof of these facts is given.
We give a complete and independent proof of both these claims in
Propositions \ref{3points} and \ref{2points}. Notice that Theorem
\ref{final} deliberately omits the case $(p,l)=(0,n+3)$. As far as
we know, the only known answer is for $n=3$ and it is given by
Wakeford in \cite{Wakeford} (see Section \ref{twisted}).

The non-existence result deserves a special comment. For $p\geq 4$ and
$l\geq 2$ Proposition \ref{donotexist} states that in $\PP^n$ no rnc
exists passing through $p$ generic points and $(n-1)$-secant to
$l$ generic codimension two linear spaces. Thus, for $p+l=n+3, p\geq 4$
and $l \geq 2$ the count of conditions always fails. The proof of this
fact is very simple: one sees that a curve with the required
properties {\it must} be reducible and hence it can not be a rnc.
At this point, one can think to allow degenerations of rnc's and
gain existence also in these cases. But this does not happen. Even
in $\PP^3$, allowing degenerations is not enough. In fact, for
$(p,l)=(4,2)$ the degree 3 curves passing through the
points and having the lines as chords split as the union of a
{\it not} intersecting conic and a line. Thus, a degree $3$ curve
of arithmetic genus $-1$. In conclusion, the non-existence is a
fact deeply related to the nature of the problem, not only coming
from a restrictive choice of curves.

We already mentioned an application of our results to the study of
higher secant varieties of some Segre-Veronese varieties. More
precisely, we consider a scheme $Z\subset\PP^n$ of $n+2$ double
points union a codimension two linear space. Using rnc's we can
easily show that $H(Z,4)$ is at least one less than expected (see
Lemma \ref{applemma}). We think that this method can be
successfully applied for studying other families of schemes
supported on points and linear spaces and we plan to investigate
this in the future. Another application involves the study of
projective equivalence of some special family of subsets and it is
described in Section \ref{projeq}. In particular, we give a
criterium to establish whether two ordered subset of $\PP^n$ each
consisting of $p$ points and $n+3-p$ codimension two linear spaces
are projectively equivalent.

For the convenience of the reader, we give an outline of the
paper. In Section \ref{definition} we give the necessary
definition and we recall some basic facts about rnc's that we will
extensively use. Also, we prove Lemma \ref{generalLEMMA} which is
the technical core of the paper.  In Section \ref{twisted} we give
a historical account of the classic results for twisted cubic
curves. Section \ref{main} contains the main results of the paper.
Finally, in Sections \ref{appsection} and  \ref{remsection} we
give some applications of our results and we make some final
remarks on the problem.

The authors wish to thank C. Ciliberto, M. Mella and F. Russo for
the interesting discussions on the topic. In particular, F.
Russo's support was crucial in starting and developing this work.
The first author had the occasion of meeting all this people at
the Workshop on Cremona Transformation held at the Politecnico di
Torino in September 2005, organized by G. Casnati, R.
Notari, and M.L. Spreafico. A special thank to the organizers for providing such a
productive scientific occasion.

\section{Notation and basic facts}\label{definition}

A {\it rational normal curve} (a {\it  rnc} for short) in $\PP^n$ is an
irreducible, reduced, smooth, rational, linearly normal curve.

A linear space $\Lambda\subset\PP^n$ is said to be {\it $(n-1)$-secant}
to a rnc $\mathcal{C}$ if $\Lambda\cap \mathcal{C}$ is a set of $n-1$
distinct points, and we will also say that the
curve $\mathcal{C}$ is {\it $(n-1)$-secant} to $\Lambda$.

Let $\mbox{Hilb}^{nt+1}(\PP^n)$ be the Hilbert scheme
parameterizing subschemes of $\PP^n$ having Hilbert polynomial
$nt+1$. Rational normal curves correspond to the points of a
smooth, irreducible, open subscheme of $\mbox{Hilb}^{nt+1}(\PP^n)$
which we denote by $\mathcal{H}$. For more on this see
\cite{MinniRagni} and the references there. We recall that
$\dim\mathcal{H}=(n-1)(n+3)$.

In this paper, we often invoke Bezout type arguments. In
particular, we are interested in showing that a given hypersurface
$X$ contains a rnc $\mathcal{C}$. This is in turn equivalent to
show that a two variable polynomial $F$ of degree $(\deg
X)\cdot(\deg\mathcal{C})$ has too many roots.

 We recall that if
$P\in \mathcal{C}$  is a multiple point for $X$, then
$F$ has a multiple root of at least the same multiplicity.
Moreover, if $X$
contains a hyperosculating space (e.g. the tangent space, the
ordinary osculating space, etc. ) to $\mathcal{C}$ in $P$, then
$F$ has again a multiple root of the proper multiplicity.

Given a variety $X$ and a natural number $n$, $(X)^n$ will denote
the product $\underbrace{X\times\ldots\times
X}_{n-\mathrm{times}}$. If $X$ is embedded in $\PP^n$ we follow
Harris in \cite[pg. 90]{Harris} and we consider the {\it
$k$-secant map}
\[
(X)^{k+1}\dashrightarrow G(k,\PP^n)
\]
mapping $k+1$ generic points of $X$ to the point of the Grassmannian $G(k,\PP^n)$
corresponding to the linear space that they
span. In particular, we denote by $\mathbb{S}^{k}X$ the closure of
the image of this map and we call it the {\it abstract variety of
secant $k$-spaces to $X$}.

We recall that a generic $2\times n$ matrix $\mathsf{M}$ of linear forms on
$\PP^n$ defines a rnc via its maximal minors.
A generalized row of $\mathsf{M}$ is any row of a matrix conjugate
to $\mathsf{M}$; similarly for a generalized column. Notice that
the zero locus of a generalized row of $\mathsf{M}$ is a point of
the related rnc $\mathcal{C}$ and viceversa. Also, the zero locus
of a generalized column intersect $\mathcal{C}$ in a zero
dimensional scheme of degree $n-1$ and also the converse holds
 as shown in the following lemma (for an alternative proof see \cite[pg.
 102]{Harris}).

\begin{lem}\label{codim2secant}
Let $\mathsf{M}$ be a $2\times n$ generic matrix of linear forms
and let $\mathcal{C}\subset\PP^n$ be the rnc defined as the rank
one locus of $\mathsf{M}$. Then, $\Lambda$ is a codimension
two linear space intersecting $\mathcal{C}$ in a degree $n-1$
scheme if and only if $\Lambda= \{F=G=0\}$ and ${F\choose G}$ is a generalized column of
$\mathsf{M}$.
\end{lem}
\begin{proof}
Let
\[
\mathsf{M}=
\left(
\begin{array}{cccc}
F_1 & F_2 & \ldots & F_n\\
G_1 & G_2 & \ldots & G_n
\end{array}
\right).
\]
If ${F\choose G}$ is a generalized column of $\mathsf{M}$, we can
substitute $\mathsf{M}$ with a conjugate matrix of the form
\[
\left(
\begin{array}{cccc}
F & F_2 & \ldots & F_n\\

G & G_2 & \ldots & G_n
\end{array}
\right).
\]
Delete the first column and consider the rank one locus of the
resulting matrix, call it $X$. Then, $X$ is a degree $n-1$ surface
and $X\cap \Lambda$=$\mathcal{C}\cap \Lambda$. Hence $\{F=G=0\}\cap
\mathcal{C}$ is a zero dimensional scheme of degree $n-1$.

Conversely, assume that $\Lambda$ is a codimension two linear
space intersecting $\mathcal{C}$ in a degree $n-1$ scheme. For
simplicity assume that $\Lambda\cap \mathcal{C}$ is a smooth set of
points, say $\{P_1,\ldots,P_{n-1}\}$. Evaluate $\mathsf{M}$ in
$P_i$ and let $V_i\subset\mathbb{C}^n$ be the space of solution of
the linear system $\mathsf{M}_{|P_i}\underline{\lambda}=0$.
Clearly, $\bigcap_{i=1}^{n-1} V_i\neq 0$ and let
$(\lambda_1,\ldots,\lambda_n)$ be a common solution. Then the generalized column
\[
{F\choose G}=\sum_1^n \lambda_i{F_i\choose G_i}
\]
is such that $\Lambda=\{F=G=0\}$.
\end{proof}

In order to prove our results, we fix some notation and we derive
a crucial technical fact playing a key role in this paper. Given
natural numbers $l$ and $p$, we consider the {\it data} space
\[
\mathcal{D}=G(n-2,\PP^n)^l\times(\PP^n)^p
\]
parameterizing sets consisting of $l$ codimension 2 linear spaces
and $p$ points in $\PP^n$. Notice that
$\dim\mathcal{D}=np+2l(n-1)$. We call an element of $\mathcal{D}$
a {\it datum}. Given a rnc $\mathcal{C}$ and a datum
$\delta=(\Lambda_1,\ldots,\Lambda_l,P_1,\ldots,P_p) \in \mathcal{D}$, we say that
$\mathcal{C}$ {\it satisfies} $\delta$ if

\[P_i\in\mathcal{C} \mbox{ for } 1\leq i\leq p\]
and

\[\Lambda_i\cap \mathcal{C} \mbox{ has degree } n-1 \mbox{ for }
1\leq i\leq l.\]

Then we consider the incidence correspondence
$\Sigma\subset\mathcal{H}\times\mathcal{D}$ defined as
\[
\Sigma=\lbrace (\mathcal{C},\delta): \mathcal{C}\mbox{ satisfies
}\delta\rbrace.
\]
With this notation, we can rephrase our problem about the
existence of rnc's: given a generic datum are there rnc's
satisfying it? If we let $\phi:\Sigma\rightarrow\mathcal{D}$ be
the natural projection, this question reduces to the following: is
$\phi$ dominant?

Finally we can introduce the main technical tool of the paper.
\begin{lem}\label{generalLEMMA}
If $p+l=n+3$ and there exists a datum $\delta\in\mathcal{D}$ such
that $\phi^{-1}(\delta)$ is a finite number of points, then $\phi$ is
dominant. Moreover, if the datum
$\delta=(\Lambda_1,\ldots,\Lambda_l,P_1,\ldots,P_p)$ is such that $\Lambda_i$ is
$(n-1)$-secant to the curve $\phi^{-1}(\delta)|_{\mathcal{H}}$ for
$1\leq i\leq l$, then the same holds for the generic element of
$\mathcal{D}$.
\end{lem}
\begin{proof}
Notation as above. First we will show that $\Sigma$ is
irreducible. Let $\psi:\Sigma\rightarrow\mathcal{H}$ be the
projection map and consider $\mathcal{C}\in\mathcal{H}$. Notice
that
\[
\psi^{-1}(\mathcal{C})_{|\mathbb{P} ^n}\simeq \mathcal{C},
\]
\[
\psi^{-1}(\mathcal{C})_{|G(n-2,\PP^n)}\simeq
\mathbb{S}^{n-2}\mathcal{C},
\]
where $\mathbb{S}^{n-2}\mathcal{C}\simeq
\mathcal{C}\times\ldots\times \mathcal{C}=(\mathcal{C})^{n-1}$.
Hence $\psi$ has irreducible fibers all having dimension
$p+(n-1)l$. Thus $\Sigma$ is irreducible and
$\dim\Sigma=p+(n-1)(l+n+3)$. Notice that, as $p+l=n+3$,
$\dim\mathcal{D}=np+2l(n-1)=\dim\Sigma$. Then $\phi$ is readily
seen to be dominant as $\mbox{Im}\phi$ is irreducible and such
that $\dim\Sigma-\dim\mbox{Im}\phi\leq\dim\phi^{-1}(\delta)=0$.

Then, let
\[
\Sigma^\circ=\lbrace(\mathcal{C},(\Lambda_1,\ldots,\Lambda_l,P_1,\ldots,P_p))\in\Sigma
: \mathcal{C}\cap \Lambda_i\mbox{ is not smooth for some }i\rbrace
\]
and notice that $\Sigma^\circ$ is closed and proper in $\Sigma$.
Hence $\dim\Sigma^\circ<\dim\Sigma$ and the second assertion
follows as the fiber of $\phi$ over the generic datum can not be
contained in $\Sigma^\circ$.
\end{proof}

\section{The $\PP^3$ case}\label{twisted}

In this section, we briefly illustrate our problem in the well
known case of rnc's in $\PP^3$. Twisted cubics have been
thoroughly investigated classically and we recall some of the many
interesting results, which are usually obtained via ad hoc
techniques.

The parameter space of twisted cubic curves has dimension $12$ and
fixing one point in $\PP^3$ imposes two conditions. In particular,
from a numerical point of view, the condition of passing through a
fixed point is equivalent to the one of touching a fixed line in
two (not fixed) points. In conclusion, given $p$ points $P_1, \ldots, P_p$ and $l$
 lines $\Lambda_1, \ldots, \Lambda_l$ in generic position in $\PP^3$, such that $p+l=6$, we expect
to find a finite
number of twisted cubic curves passing through the points and
2-secant to the lines.

For $p=6,l=0$, an answer can be obtained by considering quadric cones.
Namely, let $\mathcal{Q}$ and $\mathcal{Q}'$ be two
quadrics containing the points $P_1,\ldots ,P_6$ and with a double
point in $P_1$ and $P_2$, respectively.
The complete
intersection $\mathcal{Q} \cap \mathcal{Q}'$ is the union of the
line $P_1P_2$ and of a twisted cubic through $P_1,\ldots ,P_6$.
Hence the rnc exists and, by Bezout, it is easy to show that it is unique.

For $p=5,l=1$, we consider again quadrics.
Let $\mathcal{Q}$ and $\mathcal{Q}'$ be two smooth
quadrics containing the points and the line. Then the complete
intersection $\mathcal{Q} \cap \mathcal{Q}'$ is the union of the
given line and of the unique twisted cubic with the required properties.

For $p=4,l=2$, we
expect to find a twisted cubic passing through  four given points and
2-secant to  two given lines. But such a curve does not exist. To see
this, simply consider the unique quadric $\mathcal{Q}$ containing
the lines and  the points $P_1,P_2,P_3$; notice that, by genericity, the
fourth point $P_4$ is not on $\mathcal{Q}$. Clearly, any curve
with the required properties would be contained in $\mathcal{Q}$
by Bezout, and $P_4$ can not be a point of the curve, hence a
contradiction. Thus the naive numeric count can not be blindly
trusted any more.

Now consider the case $p=3,l=3$. Given three points and three
lines, we numerically expect to find a curve passing through the
points and having the lines as chords. But the previous situation
suggests that this might not be the case. Strangely enough, the
numerics works again and the curve exists. To see this, let
$\mathcal{Q}$ be the unique (smooth) quadric containing the  two lines
$\Lambda_1,\Lambda_2$,
and the points $P_1,P_2,P_3$. We can assume the lines to be of type $(1,0)$ on
$\mathcal{Q}$. Let $R_1,R_2$  be the points  $\mathcal{Q}\cap \Lambda_3$. The vector space of
curves of type $(1,2)$ on $\mathcal{Q}$ has dimension $6$ and
hence there exists a (unique) rnc $\mathcal{C}$ containing
$R_1, R_2, P_1,P_2, P_3$. Thus $\mathcal{C}$ is the
required twisted cubic and the numeric count works again.  It is worth of noting that this can
also be seen using the classical projective generation of the rnc
also known as Steiner's construction (see, e.g., \cite{Todd} and
\cite{Harris}).

Next, the case $p=2,l=4$. In \cite{Wakeford}, Wakeford treated
this case using a Cremona transformation. Namely, consider the
linear system of cubics containing the lines. The corresponding
map is a Cremona transformation of type $(3,3)$ mapping the
required twisted cubic curves in lines and viceversa. Existence
and uniqueness follow simply by taking the preimage of a line.
Notice that this is again a constructive method: let $\mathcal{S}$
and $\mathcal{S}'$ be the cubics containing the lines and one of
the points. Then $\mathcal{S}\cap \mathcal{S}'$ splits as the
union of the four lines, the two four secant to them and a
residual twisted cubics with the required properties.

Finally the $p=0,l=6$ case. This has been studied again by
Wakeford in \cite{Wakeford} via a Cremona transformation and a
chords argument. First observe that two twisted cubic can have at most ten common
chords and this is the case if they are generic. To see this use
the linear system of quadrics through the first curve to map
rationally $\PP^3$ onto $\PP^2$. This map contracts all the chords
of the first curve to points. The second curve maps to a degree
six rational curve whose double points corresponds to common
chords. Hence, ten common chords exist. Now apply the Cremona
given by the linear system of cubic surfaces containing the four generic
lines $\Lambda_1, \ldots, \Lambda_4$: the two extra lines $\Lambda_5,  \Lambda_6$ are mapped to twisted cubic
curves having four common chords. The preimages of the remaining six chords give
six twisted cubic curves 2-secant to the six lines  $\Lambda_1, \ldots, \Lambda_6$. Notice
that this is the only case in which we have existence but {\it
not} uniqueness. The existence of more than one curve makes the
problem considerably harder in higher dimension and in fact it
still remains unsolved.

Using quite ad hoc and special
arguments, the considerations above give a complete description of the
situation in $\PP^3$, which we summarize in the following

\begin{prop}\label{P3}
 In $\PP^3$ consider  $p$ points and $l$ lines in generic position
such that $p+l=6$. Then there exists a rational normal curve
passing through the points and 2-secant to the lines for
\[
(p,l)=(6,0),(5,1),(3,3),(2,4),(1,5),(0,6).
\]
In the case
$(p,l)=(4,2)$ the curve does not exist. Moreover, the curve is
unique in all cases but the $(p,l)=(0,6)$ case, where six such
curves exist.
\end{prop}

\section{General results}\label{main}

In this section we extend to  $\PP^n$   the results  of Proposition \ref{P3}
for every $(p,l)$ such that $p+l=n+3$, $n >3$, and $p \geq 1$. The case
$(p,l) = (0,n+3)$ is still open.

The results of this section  follow the
paradigm ``given $p$ points and $l$ lines in {\it generic}
position" then  ``some conclusions follows". We mainly use
Lemma \ref{generalLEMMA} and we  show that a Zariski non-empty open subset of
the appropriate data space $\mathcal{D}$ exists, and the proper
conclusion holds for all data in that subset.

\begin{thm}\label{castelnuovo} Given  $n+3$ points in $\PP^n$ in generic position,
there exists a unique rational normal curve passing through them.
\end{thm}

\begin{proof} A constructive proof can be found in \cite[pg. 10]{Harris} and
\cite{Bordiga}. Here we give a classic proof via Cremona
transformations.

Let $P_1, \ldots,P_{n+3}$ be the given points.
We may assume that $P_1, \ldots,P_{n+1}$  are the
coordinate points. Consider the linear system of degree $n$ hypersurfaces having the
coordinate points as singular points of multiplicity $n-1$. If we
denote the map associated to the linear system as $\varphi$, it is
well known that it is a Cremona transformation of type $(n,n)$. In
particular, $\varphi$ maps rnc's through the coordinate points in
lines and viceversa. Hence, the preimage of the unique line
joining $\varphi(P_{n+2})$ and $\varphi(P_{n+3})$ is the required rnc.

\end{proof}

\begin{prop}\label{onePOINTlessPROP}
Consider  $n+2$ points in $\PP^n$, and a codimension two linear
space in generic position. Then, there exists a unique rational
normal curve passing through the points and $(n-1)$-secant to the
linear space.
\end{prop}
\begin{proof}

We use Lemma \ref{generalLEMMA} and its notation. Thus, to show
existence we have to produce a datum $\delta$ such that
$\phi^{-1}(\delta)$ is a single point.

Let
\[
\mathsf{M}= \left(
\begin{array}{ccc}
F_1 & \ldots & F_n \\
G_1 & \ldots & G_n
\end{array}
\right)
\]
be a generic $2\times n$ matrix of linear forms
and denote by $\mathcal{C}$ the rnc defined by its $2\times 2$
minors.

By the genericity of $\mathsf{M}$, we have that
\begin{itemize}

\item for a generic choice of pairs
$(a_j,b_j)\in\mathbb{C}^2,j=1,\ldots,n+2$, the points $P_j=\{a_j F_1+b_j
G_1=\ldots =a_j F_n+b_j G_n=0\}$ are distinct;

\item $\Lambda=\{F_1=G_1=0\}$ is a codimension two linear space.
\end{itemize}

Now consider the datum
$\delta=(\Lambda,P_1,\ldots,P_{n+2})\in\mathcal{D}$. Clearly, $\delta$
is in the image of $\phi$. We will now show that
$\phi^{-1}(\delta)$ consist of a single
point. Let $\mathcal{Q}_{j}$ be the quadric defined by $\left|\begin{array}{cc}F_1 & F_j \\
G_1 & G_j\end{array}\right|=0,j=2,\ldots,n$. Observe
that $\Lambda \subset \mathcal{Q}_{j}$ and $P_1,\ldots,P_{n+2}\in
\mathcal{Q}_{j}$. Moreover, a simple rank argument yields

\[\bigcap_j \mathcal{Q}_{j}=\Lambda\cup \mathcal{C}.\]

In fact, if a point $P$ lies in this intersection but not in $\Lambda$,
then all the columns of $\mathsf{M}$ evaluated in $P$ are
proportional to the first column and hence $P\in \mathcal{C}$.
Now it is easy to check that any rnc $\mathcal{C}'$ satisfying the datum $\delta$ is
contained in $\mathcal{Q}_{j},j=2,\ldots,n$ by Bezout, and hence
$\mathcal{C}'$ and $\mathcal{C}$ coincide. In conclusion,
$\phi^{-1}(\delta)_{|\mathcal{H}}=\mathcal{C}$ and the map $\phi$
is dominant.

To prove uniqueness, it suffices to repeat the Bezout argument
above.
\end{proof}

\begin{prop}\label{donotexist}
Let $p$ and $l$ be integers such that $p\geq 4$ and $l\geq 2$.
Consider  $p$ points and $l$ codimension two linear
spaces in generic position in $\PP^n$. Then
 no rational normal curve passing through the points and
$(n-1)$-secant to the linear spaces exists.
\end{prop}
\begin{proof}
It is enough to prove the statement for $p=4$ and $l=2$. Let $\Lambda_1$
and $\Lambda_2$ be the linear spaces and let $P_1,\ldots,P_4$ be the
points. We want to show that there exists a quadric containing the
scheme $X=\Lambda_1\cup \Lambda_2\cup P_1\cup P_2\cup P_3$, i.e. we want to
show that $h^0 \mathcal{I}_X(2)>0$. The linear space $\Lambda_1$ imposes
${n\choose 2}$ independent conditions on quadrics and notice that
$\Lambda_1\cap \Lambda_2\simeq\PP^{n-4}$. Thus $\Lambda_1\cup \Lambda_2$ imposes
$2{n\choose 2}-{n-2\choose 2}$ conditions. In conclusion
\[
h^0\mathcal{I}_X(2)\geq h^0\mathcal{O}_{\PP^n}(2)-\left[2{n\choose
2}-{n-2\choose 2}\right]-3=1.
\]

Let $\mathcal{Q}$ be a quadric containing $X$ and notice that
$P_4\not\in \mathcal{Q}$ by genericity. Suppose that a rational
normal curve with the requires properties exists, say
$\mathcal{C}$. We will show that $\mathcal{Q}\supset \mathcal{C}$,
hence a contradiction. Let $t$ be the degree of the scheme
$\mathcal{C}\cap \Lambda_1\cap \Lambda_2$ and notice that $\mathcal{Q}$ is
singular along the intersection $\Lambda_1\cap \Lambda_2$. Hence the
degree of $\mathcal{Q}\cap \mathcal{C}$ is at least
\[
3+(n-1-t)+(n-1-t)+2t=2n+1
\]
so by Bezout we get $\mathcal{Q}\supset
\mathcal{C}$, a  contradiction.
\end{proof}

\begin{prop}\label{3points}
Consider in $\PP^n$ three points and $n$ codimension two linear
spaces in generic position. Then, there exists a unique rational
normal curve passing through the points and $(n-1)$-secant to the
linear spaces.
\end{prop}
\begin{proof}
Notation as in Lemma \ref{generalLEMMA}.
Existence is proved if we produce a datum $\delta$ such that
$\phi^{-1}(\delta)$ is a single point.

Let  $\mathsf{M}$ and $\mathcal{C}$ be
as in the proof of Proposition \ref{onePOINTlessPROP}.

By the genericity of $\mathsf{M}$, we have that
\begin{itemize}
\item $\{F_1=\ldots=F_n=0\}$ is a point, say $P_1$;

\item $\{G_1=\ldots=G_n=0\}$ is a point, say $P_2$;

\item $\{F_1+G_1=\ldots=F_n+G_n=0\}$ is a point, say $P_3$;

\item $\{F_i=G_i=0\}$ is a codimension two linear space, say $\Lambda_i$,
$i=1,\ldots,n$.

\end{itemize}

Now consider the datum
$\delta=(\Lambda_1,\ldots,\Lambda_n,P_1,P_2,P_3)\in\mathcal{D}$. Clearly,
$\delta$ is in the image of $\phi$. We will now show that
$\phi^{-1}(\delta)$ consist of a single point. For $1\leq i<j\leq
n$, let
$\mathcal{Q}_{ij}$ be the quadric defined by $\left|\begin{array}{cc}F_i & F_j \\
G_i & G_j\end{array}\right|=0$  and notice that
$ \Lambda_i,\Lambda_j \subset \mathcal{Q}_{i,j}$ and $P_1,\ldots,P_3\in
\mathcal{Q}_{i,j}$. It easily follows that  any rnc $\mathcal{C}'$ satisfying the
datum $\delta$ is contained in all the $\mathcal{Q}_{ij}$'s by
Bezout. Hence the curves $\mathcal{C}$ and $\mathcal{C}'$ coincide
as they have the same defining ideal. Hence
$\phi^{-1}(\delta)_{|\mathcal{H}}=\mathcal{C}$, thus the map $\phi$
is dominant. This proves the  existence of the rnc with the
desired properties.  Uniqueness follows by the Bezout
argument above.
\end{proof}

\begin{prop}\label{2points}
Consider in $\PP^n$ two points and $n+1$ codimension two linear
spaces in generic position. Then, there exists a unique rational
normal curve passing through the points and $(n-1)$-secant to the
linear spaces.
\end{prop}
\begin{proof} First we will recall the classical proof, only given in $\PP^4$, and then we
produce ours.

\vskip .5cm \noindent({\it Veneroni-Todd})

The idea is to use a Cremona transformation, also known as
Veneroni's transformation. We refer to the classical work \cite{Todd} for more
details. Consider the linear system of quartics containing five
given planes in generic position and let
$\varphi:\PP^4\dashrightarrow\PP^4$ be the corresponding map. This
map can be shown to be a Cremona of type $(4,4)$ and it maps each
rnc 3-secant to each of the planes in a line and vice versa. Given
the points $P_1$ and $P_2$, existence and uniqueness follow by
considering the preimage of the unique line joining $\varphi(P_1)$
and $\varphi(P_2)$.

\vskip .5cm \noindent({\it Complete proof})

We use Lemma \ref{generalLEMMA} and its notation. Thus, the
existence part of the proof is completed if we exhibit a datum
$\delta$ such that $\phi^{-1}(\delta)$ is non-empty and finite.

Let  $\mathsf{M}$ and $\mathcal{C}$ be
as in the proof of Proposition \ref{onePOINTlessPROP}.

We have that
\begin{itemize}
\item $\{F_1=\ldots=F_n=0\}$ is a point, say $P_1$;

\item $\{G_1=\ldots=G_n=0\}$ is a point, say $P_2$;

\item $\{ F_i=G_i=0\}$ and $\{\sum F_i=\sum
G_i=0\}$ are  codimension two linear spaces, say $\Lambda_i$
(for $i=1,\ldots,n$)
and $\Lambda$, respectively.
\end{itemize}

Now consider the datum
$\delta=(\Lambda_1,\ldots,\Lambda_n,\Lambda,P_1,P_2)\in\mathcal{D}$. Clearly,
$\delta$ is in the image of $\phi$. We will now show that
$\phi^{-1}(\delta)$ consist of a single point. Let $\mathcal{C}'$
be another rnc satisfying the datum $\delta$. By Lemma
\ref{codim2secant} we know that a defining matrix of
$\mathcal{C}'$ can be chosen of the form
\[
\mathsf{M}'= \left(
\begin{array}{ccc}
a _1 F_1+b_1 G_1 & \ldots & a _n F_n+b_n G_n \\
a_1' F_1+b_1' G_1 & \ldots & a _n' F_n+b_n' G_n
\end{array}
\right).
\]

By the genericity of $\mathsf{M}$, $P_1$ annihilates all the $F_i$'s
but none of the $G_i$'s. Hence the vectors $(b_1, \ldots, b_n)$ and $(b'_1, \ldots, b'_n)$
are proportional. Thus a linear combination of the rows of
$\mathsf{M'}$ eliminates all the $G_i$'s in the first row;
similarly for the $F_i$'s in the second row using $P_2$. So that
$\mathsf{M}'$ is conjugate to a matrix of form
\[
\left(
\begin{array}{ccc}
F_1 & \ldots & F_n \\
c_1 G_1 & \ldots & c_n G_n
\end{array}
\right).
\]

Now recall that $\Lambda=\{\sum F_i=\sum G_i=0\}$ intersects both
$\mathcal{C}$ and $\mathcal{C}'$ in a degree $n-1$ scheme. Then,
by Lemma \ref{codim2secant}, the vector space $\langle\sum
F_i,\sum G_i\rangle$ contains both $\sum e_i F_i$ and $\sum e_ic_i
G_i$ for some choice of the $e_i$'s in $\CC$. Since
$F_1,\ldots,F_n$ and $\sum G_i$ are linearly independent by the
genericity of $\mathsf{M}$, we have $e_1=\ldots =e_n$.
Similarly by the independence of $G_1,\ldots,G_n$ and
$\sum F_i$, we get $e_1c_1=\ldots =e_nc_n$. These conditions force
$\mathsf{M}$ and $\mathsf{M}'$
to be conjugate, hence $\mathcal{C}$ and $\mathcal{C}'$ coincide,
 $\phi^{-1}(\delta)_{|\mathcal{H}}=\mathcal{C}$, and
the map $\phi$ is dominant.

To prove  uniqueness, it is enough to repeat the argument
above about the defining matrices.

\end{proof}

\begin{prop} \label{onepoint}
Consider in $\PP^n$ a point and $n+2$ codimension two linear
spaces in generic position. Then, there exists a
unique rational normal curve passing through the points and
$(n-1)$-secant to the linear spaces.
\end{prop}
\begin{proof}
\vskip .5cm \noindent({\it Todd})

We only sketch the original proof given in $\PP^4$ and we refer to
\cite{Todd} for more details. Let $\Lambda_1,\ldots,\Lambda_6$ be the planes
and $P$ be the point. Consider the linear system of quartics
through $\Lambda_1,\ldots,\Lambda_5$ and let $\varphi$ be the corresponding
rational map. Notice that $\varphi$ is a Cremona of type $(4,4)$.
The rnc's satisfying the data are among the preimages of the lines
through $\varphi(P)$. Since $\varphi$ maps $\Lambda_6$ in a Bordiga
surface, the unique trisecant line through $\varphi(P)$ gives the
desired curve.

\vskip .5cm \noindent({\it Complete proof})

Notation as in Lemma \ref{generalLEMMA}. As  usual,
let  $\mathsf{M}$ and $\mathcal{C}$ be
as in the proof of Proposition \ref{onePOINTlessPROP}.
We have the following:
\begin{itemize}

\item $\{F_1=\ldots=F_n=0\}$ is a point, say $P$;

\item $\{ F_i=G_i=0\}$ are codimension
two linear spaces, say $\Lambda_i$,
 $i=1,\ldots,n$.
\end{itemize}

Then consider a rnc $\mathcal{C}'$, defined by a matrix
$\mathsf{M}'$, and impose that $P\in \mathcal{C}'$ and that
$\Lambda_1,\ldots,\Lambda_n$ intersect also $\mathcal{C}'$ in a degree $n-1$
scheme.
Since $P$ annihilates all the $F_i$'s, but none of the $G_i$'s,
by arguments similar to the ones used in the proof of the previous proposition,
we get
\[
\mathsf{M}'= \left(
\begin{array}{cccc}
F_1 & F_2 & \ldots & F_n \\
G_1 & a _2 F_2+b_2 G_2 &\ldots & a _n F_n+b_n G_n
\end{array}
\right).
\]

Now we choose two extra common secant spaces to
$\mathcal{C}$ and $\mathcal{C}'$ in such a way that the curves are
forced to coincide. The proof is different depending on the parity of  $n$.

If  $n$ is odd, let $n=2m-1$. Consider the linear spaces $
\Lambda_{n+1}=\{\sum_{i=1}^m F_i=\sum_{i=1}^m G_i=0\}$ and $\Lambda_{n+2}
=\{\sum_{i=m}^{2m-1} F_i=\sum_{i=m}^{2m-1} G_i=0\}$, and require that $\mathcal{C}'$
is $(n-1)$-secant to them.  Since $\Lambda_{n+1}$ intersect
$\mathcal{C}'$ in a degree $n-1$ scheme, by Lemma
\ref{codim2secant} there exist constants
 $e_1, \ldots, e_n \in \CC $ such that
\[
\Lambda_{n+1}=\left \{\sum_{i=1}^n e_i F_i=e_1G_1+\sum_{i=2}^n e_i(a_i
F_i+b_i G_i)=0\right \}
\]

Since the vector space  $\langle\sum _{i=1}^m
F_i,\sum _{i=1}^mG_i\rangle$ contains   $\sum_{i=1}^n e_i F_i$
and the $n+1$ linear forms  $F_1, \ldots, F_n, \sum _{i=1}^mG_i$
are linearly independent,  we get $e_1= \ldots= e_m  $ and
$e_{m+1}= \ldots= e_n=0 $.
Moreover also $e_1G_1+\sum_{i=2}^n e_i(a_i F_i+b_i G_i)$ is an element of the vector space
above, then
there exist  $a$, $b \in \CC$ such that
\[
a\sum_{i=1}^m F_i + b\sum_{i=1}^m G_i=G_1+\sum_{i=2}^m (a_i
F_i+b_i G_i).
\]

This equality involves the $2m$ linear forms
$F_1, \ldots, F_m,G_1, \ldots, G_m$, which are
linearly independent. By comparing their
coefficients we get
$a=a_2=\ldots =a_m =0$ , and $b=b_2=\ldots =b_m =1$.
Analogously,  imposing that $ \Lambda_{n+2}$ is a
$(n-1)$-secant space to $\mathcal{C}'$, we get $a_i=0$ and $b_i=1$
for all $i$. Thus $\mathsf{M}=\mathsf{M}'$, the curves
coincide, and if we let $\delta=(\Lambda_1,\ldots,\Lambda_{n+2},P)$, then $\phi^{-1}(\delta)$ consists of a
single point.

In the case $n$ even, let $n=2m-2$.Now  consider the linear spaces
$\Lambda_{n+1}=\{\sum_{i=1}^m F_i=\sum_{i=1}^m G_i=0\}$ and $
\Lambda_{n+2}=\{F_1+\sum_{i=m}^{2m-2} F_i=G_1+\sum_{i=m}^{2m-2} G_i=0\}$,
and analogously to the  $n$ odd case,  require that $\mathcal{C}'$
is $(n-1)$-secant to them.  Arguing as above,
we get that there exist  $a$, $b \in \CC$ such that
\[
a\sum_{i=1}^m F_i + b\sum_{i=1}^m G_i=G_1+\sum_{i=2}^m (a_i
F_i+b_i G_i).
\]

Now we may assume that
$G_1=F_1+\ldots+F_m+G_2+\ldots+G_m$.
By comparing the coefficients of the independent linear
 forms $F_1, \ldots, F_m, G_2, \ldots, G_m$  we easily get
\[  a_2=\ldots =a_m=0 \mbox{ and }   b_2=\ldots =b_m.\]

Since  $\Lambda_{n+2}$ intersect the
curve $\mathcal{C}'$ in a degree $n-1$ scheme,
 the vector space  $\langle
F_1+\sum_{i=m}^{2m-2} F_i,G_1+\sum_{i=m}^{2m-2} G_i
\rangle$ contains   $\sum_{i=1}^n e'_i F_i$ and
$e'_1G_1+\sum_{i=2}^n e'_i(a_i F_i+b_i G_i)$
for some $e'_i$'s $\in \CC$. Repeating the usual arguments
we get that there exist
$a'$, $b' \in \CC$ such that

\begin{eqnarray}
& \nonumber
a'(F_1+F_m+\ldots+F_{2m-2})+b'(G_1+G_m+\ldots+G_{2m-2}) &
\\
&\nonumber =  G_1+b_m
G_m+\sum_{i=m+1}^{2m-2} (a_i F_i+b_i G_i) \hfill &
\end{eqnarray}
Since
\[
\{G_1,G_m,\ldots,G_{2m-2},F_{m+1},\ldots,F_{2m-2},F_1+F_m\}
\]
is a  set of $n+1$ linearly independent  forms, again by comparing their
coefficients we get
$b_m= \ldots =b_{2m-2}=1$ and  $a_{m+1}= \ldots =a_{2m-2}=0$.
Summing up these relations with the
previous ones, we obtain $\mathsf{M}=\mathsf{M}'$
and if we let $\delta=(\Lambda_1,\ldots,\Lambda_{n+2},P)$, then $\phi^{-1}(\delta)$
 consists of a single
point.

The  existence part of the proof is now
completed. To show { uniqueness} it is enough to repeat argue
as above on the defining matrices.
\end{proof}

We summarize the results of  Theorem \ref{castelnuovo},
Proposition \ref{onePOINTlessPROP} to \ref{onepoint} in
the following
\begin{thm}\label{final}
Let $n,p$ and $l$ be positive
integers such that
$$n\geq
3, \ \ p\geq 1\mbox{ and} \ \ p+l=n+3 .$$
 Choose  $p$ points in $\PP^n$ and $l$
codimension two linear spaces in generic position. Then, only for
the values
\[
(p,l)=(n+3,0),(n+2,1),(3,n),(2,n+1),(1,n+2)
\]
does there
exist a unique rational normal curve passing through the points and
$(n-1)$-secant to the linear spaces. In the other cases, that is for
 $p\geq 4$ and $l\geq 2$, no such curve exists.

\end{thm}

\section{Applications}\label{appsection}

\subsection{Postulation of schemes and defectivity}\label{postulationAPP}

Theorem \ref{final} can be used to produce
schemes that impose less conditions than expected to forms of some degree.
Here we only give an example to show the main ideas. A thorough  study will be the subject of a
forthcoming paper.

Let $P_1,\ldots,P_{n+2}$ be $n+2$ generic points in $\PP^n$ $(n>2),$
and let $\Lambda$ be a generic codimension two linear space, with defining ideals
 $I_{P_1}, \ldots, I_{P_{n+2}}, I_\Lambda$, respectively.
Consider the scheme $X$ having ideal
\[I_X=(I_{P_1})^2\cap \ldots\cap (I_{P_{n+2}})^2\cap (I_\Lambda)^2.\]
It is easy to compute the expected Hilbert function of $X$ in degree $4$, namely
$$h=(n+2)(n+1)+{n+2\choose 4}+2{n+1\choose 3}.$$
The following lemma shows that the scheme $X$ have not the expected postulation:
\begin{lem}\label{applemma}
Notation as above, the scheme $X$ does not impose the
expected number of conditions on degree 4 hypersurfaces, i.e.
$H(X,4)\leq h-1$.
\end{lem}
\begin{proof}
Consider the scheme $X'$ with defining ideal
\[I_{X'}=I_{P_1}\cap (I_{P_2})^2\ldots\cap (I_{P_{n+2}})^2\cap (I_\Lambda)^2.\]
By Theorem \ref{final} we know that there exists a rnc $\mathcal{C}$
through the $P_i$'s and having $\Lambda$ as a $(n-1)$-secant space.
Moreover, any
element $F\in (I_{X'})_4$ vanishes on $\mathcal{C}$ by a standard
Bezout argument as the degree of $\{F=0\}\cap \mathcal{C}$ is
\[ 1+2(n+1)+2(n-1)>4n.\]
Hence, all quartic hypersurfaces through $X'$ have a fixed tangent direction in
$P_1$. This is enough to conclude that $H(X,4)$ is at least one
less than expected.
\end{proof}
As a straightforward application of Lemma \ref{applemma} the next corollary
gives a not trivial
defectiveness result for Segre-Veronese varieties. The same statement
can be deduced by the classification given in \cite{AAdlandsvik},
where {\AA}dlandsvik  heavily uses his theory of the
joint of varieties.
Abrescia in \cite{Abrescia} proposes a
simplified proof of this result.
\begin{cor}\label{applcor} Let V be the Segre-Veronese variety $\PP^1\times\PP^n$ embedded
with bidegree $(2,2)$. Then $V$ is $(n+1)$-defective, i.e. the
secant variety $S^{n+1}(V)$ has not the expected dimension.
\end{cor}
\begin{proof}
The conclusion follows from the previous lemma and the results of
Section 1 in \cite{CGG3}, where the authors relates the
$(n+1)$-defectiveness of $V$ with the Hilbert function in degree
$4$ of schemes consisting of $n+2$ double points and a double
codimension two linear space.
\end{proof}

\subsection{Projectively equivalent subsets}\label{projeq}

It is well known that any two ordered sets of $n+2$ generic points in
$\PP^n$ are projectively equivalent, i.e. there is an  automorphism
of $\PP^n$ mapping one set in the other preserving the order
of the points.

If  $\mathbb{X},\mathbb{Y}$ are two ordered sets of $n+3$ points in $\PP^n$,
 it is
interesting to look for conditions assuring the projective
equivalence. Theorem \ref{castelnuovo} gives the answer:
via the unique rational normal curves passing through the
points of each set, we  map $\mathbb{X}$ and $\mathbb{Y}$ into
$\PP^1$. Then the question is answered via cross ratios.
For more than  $n+3$ points, the problem stays open (see \cite[pg.
8]{Harris}).

If linear spaces and points are taken into account,
as far as we know, there are no similar results. Our Theorem
\ref{final} allows to give an answer for a special family of
subsets of $\PP^n$. Namely, the sets consisting of $p$ points and
$l$ codimension two linear spaces in generic
position
with $p+l=n+3, n \geq 3, p\geq 1$.

For example, consider in $\PP^3$ a generic set
$A$ consisting of three lines and
three points.
Via the unique rnc through the three points and 2-secant to the three lines,
 we obtain  a subset
$A'$  of nine points in $\PP^1$. Then, another set
 $B$ of three lines and
three points is projectively equivalent to $A$ if and only if $A'$ and
$B'$ are so, where $B'$ is similarly constructed. Thus, an answer can be obtained again via cross
ratios.

\section{Final remarks}\label{remsection}

\subsection{The case of $n+3$ codimension two linear spaces}

In this situation we do not have a complete answer to the basic
question: given in $\PP^n$ $n+3$ linear spaces of codimension two
in generic position, are there rational normal curves
$(n-1)$-secant to these
 spaces? For $n=3$, the answer is
positive (see Proposition \ref {P3}), and we  also have a proof
 using Lemma \ref{generalLEMMA} and a cubic surface having five
of the six lines as exceptional divisors.
Unfortunately this proof does not extend to
$n>3$.

\subsection{Mixed conditions}

In this paper we generalized the classical Theorem
\ref{castelnuovo} by substituting points with codimension two
linear spaces, motivated by a count of conditions. If we look for
further generalizations, we can again rely on a count of
conditions for inspiration. For example, consider in $\PP^{2m+3}$
$m+4$ points, a linear space of dimension $m$, and a linear space
of dimension $m+1$. Then, by counting conditions, we expect to
find a rnc passing through the points, $(m+1)$-secant to the $m$
dimensional space and $(m+2)$-secant to the $m+1$ dimensional
space. Actually, we can prove that such a curve exists and it is
unique. Note that for $m=1$ this result yields a statement similar
to the one of Lemma \ref{applemma}. By this result, analogously to
Corollary  \ref{applcor}, we deduce another proof for the
classically known 5-defect of the Segre-Veronese
$\PP^2\times\PP^3$ embedded with bidegree $(1,2)$ (for a modern
proof see for instance \cite{CaCh}, Theorem 4.3). We are presently
studying in  which direction we have to move in order to get a
result as comprehensive as Theorem \ref{final}.

\subsection{The higher dimensional case}

The title of the present paper can be rephrased as: existence
results for 1-dimensional Veronese varieties. Thus, it is
extremely natural to pose questions similar to the ones we
addressed here in the higher dimensional case. In dimension 2,
there is a well known result by Kapranov \cite{Kapranov} about the
existence of Veronese surfaces containing special sets of points.
Recently, Graber and Ranestad re-proposed Kapranov's result and
improved it by considering the existence of Veronese surfaces
``well intersecting" a special configuration of linear spaces.
These are only partial results and the problem remains open even
in the two dimensional case. As far as we know, no relevant
results exist in higher dimension.

\bibliographystyle{alpha}
\bibliography{carlini}

\end{document}